\documentclass[12pt]{article}
\usepackage{latexsym}

\def\bs{\begin{subequations}}
\def\es{\end{subequations}}

\catcode`\@=11

\newtoks\@stequation
\def\subequations{\refstepcounter{equation}
  \edef\@savedequation{\the\c@equation}%
  \@stequation=\expandafter{\theequation}%   %only want \theequation
  \edef\@savedtheequation{\the\@stequation}% % expanded once
  \edef\oldtheequation{\theequation}%
  \setcounter{equation}{0}%
  \def\theequation{\oldtheequation\alph{equation}}}

\def\endsubequations{\setcounter{equation}{\@savedequation}%
  \@stequation=\expandafter{\@savedtheequation}%
  \edef\theequation{\the\@stequation}\global\@ignoretrue}

\makeatletter%  allow access to internal LaTeX commands
        \renewcommand{\theequation}{\thesection.\arabic{equation}}%
        \@addtoreset{equation}{section}%
\makeatother%  turn off access to internal LaTeX commands

\renewcommand{\thefootnote}{\fnsymbol{footnote}}

\begin{document}

\begin{titlepage}

August 18, 2008 

Added another proof in Section 2 and also a new Appendix

\begin{center}        \hfill   \\
            \hfill     \\
                                \hfill   \\

\vskip .25in

{\large \bf Integration with Functions of a Quaternionic Variable \\}

\vskip 0.3in

Charles Schwartz\footnote{E-mail: schwartz@physics.berkeley.edu}

\vskip 0.15in

{\em Department of Physics,
     University of California\\
     Berkeley, California 94720}
        
\end{center}

\vskip .3in

\vfill

\begin{abstract}

Recent innovations in the differential calculus for functions 
of non-commuting variables, beginning with a quaternionic variable,
are now extended to consider some integration. 

\end{abstract}

\vfill

\end{titlepage}

\renewcommand{\thefootnote}{\arabic{footnote}}
\setcounter{footnote}{0}
\renewcommand{\thepage}{\arabic{page}}
\setcounter{page}{1}

\section{The Differential}%1

In a recent paper \cite{1} I showed how to expand 
\begin{equation}
F(x+\delta) = F(x) + {\cal{D}}F(x) + O(\delta^{2}) \label{a0}
\end{equation}
when both $x$ and $\delta$ were general quaternionic variables, and 
thus 
did not commute with each other; and we found the general formula,  
\begin{equation}
 {\cal{D}} F(x) = F^{\prime}(x) \delta_{\parallel} + [F(x) - F(x^{*})] 
(x-x^{*})^{-1}\;\delta _{\perp}.\label{a1}
\end{equation}
The specifics of how to construct the two components of 
$\delta=\delta_{\parallel}+\delta_{\perp}$ are, 
\begin{equation}
\delta_{\parallel} = \frac{1}{2}(\delta - u_{x}\;\delta\; u_{x}), 
\;\;\;\;\;\delta_{\perp} = \frac{1}{2}(\delta + u_{x}\;\delta\; u_{x}),
\label{a2}
\end{equation}
where $u_{x}$ is a unit imaginary that depends upon the location of 
$x = \xi_{0} + i\xi_{1} + j\xi_{2} + k\xi_{3}$  as a point 
in a four-dimensional Euclidean space: 
\begin{equation}
u_{x} = (i\xi_{1} + j\xi_{2} + k\xi_{3})/r, \;\;\;\;\; r = \sqrt{\xi_{1}^{2} 
+ \xi_{2}^{2}+\xi_{3}^{2}}. \label{a3}
\end{equation}

In another recent paper \cite{2} similar results were obtained for  
a different kind of 
non-commuting variables - matrices over the complex numbers.

This may be called the beginning of non-commutative calculus - the 
differential part.  Now we want to try looking at the integral part 
of that calculus.

\section{The Integral}%2

In ordinary calculus of functions of the real variable t, we know 
what is meant by an integral, such as $\int f(t)\;dt$. But when we 
consider quaternionic (or other non-commuting) variables it is 
unclear even how to write such an expression. (See further, Appendix 
A.)

Alternatively, we can start with the 
defining relation between the integral and the differential:
\begin{equation}
 \int_{a}^{b}\; df(t) = \int_{a}^{b} \frac{df(t)}{dt}\;dt=f(b)-f(a) ;\label{b1}
\end{equation}
and this is what we shall generalize for our non-commuting quaternionic 
variable $x$ as, 
\begin{equation}
\int_{a}^{b}\; {\cal{D}} F(x) = F (x_{b}) - F(x_{a}).\label{b2}
\end{equation}

We  define this integral as an additive operation along a path in 
that four-dimensional space of the real variables $\xi$, 
\begin{equation}
x = x_{path}(s), \;\;\;\;\; x_{path}(0) = x_{a}, \;\;\;\;\; 
x_{path}(1) = x_{b}
\end{equation}
where $s$ is a real continuous parameter.

Next, we  subdivide 
that path, whatever it may be, into a large number of infinitesimal 
increments. 
\begin{equation}
\int_{a}^{b} = \sum_{n=1}^{n=N} \;\int^{(n)}, \;\;\;\;\;\;\;\;   
\int^{(n) }=\int_{x_{n-1}}^{x_{n}}, \;\;\;\;\; n=1, 
\ldots,N\label{b3}
\end{equation}
where $x_{0} = x_{a}$ and $x_{N} = x_{b}$.

In any segment of this path we choose the line of integration, with 
the integrand ${\cal{D}} F(x)$,  to be 
the sum of two infinitesimal parts: 
\begin{equation}
x_{n} - x_{n-1} = \delta = \delta_{\parallel} + \delta_{\perp}.
\end{equation}
The first part is  ``parallel'' to 
the direction of $x$ at that point, giving the contribution
\begin{equation}
\int_{\parallel}\;{\cal{D}} F(x) = F^{\prime}(x) \delta_{\parallel}.\label{b4}
\end{equation}
Then the second part is ``perpendicular'', giving the contribution
\begin{equation}
\int_{\perp}\; {\cal{D}} F(x) = [F(x)-F(x^{*})](x-x^{*})^{-1}\; 
\delta_{\perp}.\label{b5}
\end{equation}

The sum of these two parts is thus nothing other than
\begin{equation}
F(x_{n}) - F(x_{n-1})\label{b6}
\end{equation}
to first order in the interval $\delta$.  The entire sum  then 
  results in Eq. (\ref{b2}).
 
As an alternative to this ``staircase'' construction of the 
integration, we can look at the following model.  Say that the path 
of integration is a straight line
\begin{equation}
x = \alpha+i\beta + j\gamma s \label{b7}
\end{equation}
where $s$ is a real variable that goes from $0$ to $1$.  Then we use 
the definitions Eqs. (\ref{a2}), (\ref{a3}) to calculate the 
quantities in ${\cal{D}}F(x)$.  This is a rather tedious procedure; 
but I have carried it out for $F(x) = x, x^{2}, x^{3}$ and found that the 
formula Eq. (\ref{b2}) is verified. 

Another general proof can proceed as follows. If we start with the 
coordinate along the path $x(s) = \xi_{0}(s)+i\xi_{1}(s) + 
j\xi_{2}(s) + k\xi_{3}(s)$, then we can simply write,
\begin{equation}
{\cal{D}} x(s) = ds \frac{dx(s)}{ds}
\end{equation}
since there is no commutativity problem in this representation. It is 
also true that we can express any function as 
\begin{equation}
F(x(s)) = A(s) + B(s)\;x(s)
\end{equation}
where $A$ and $B$ are real functions, the only quaternions being in 
the single factor $x(s)$.  We then see that the integral becomes 
quite ordinary:
\begin{equation}
\int_{a}^{b}\; {\cal{D}} F(x(s)) = \int_{0}^{1}\; ds 
\frac{dF(x(s))}{ds} = F(x(s))|_{0}^{1} = F(x_{b}) - F(x_{a}).
\end{equation}

In \cite{1} it was shown that this differential operator ${\cal{D}}$ 
obeys the Leibnitz rule; and thus we get the identity, usually 
called ``integration by parts'',
\begin{equation}
\int_{a}^{b} F(x)\;{\cal{D}} G(x) = F(x_{b})G(x_{b}) - 
F(x_{a})G(x_{a}) - \int_{a}^{b}\;({\cal{D}}F(x))\;G(x). \label{b8}
\end{equation}

Loosly speaking, integration is the inverse of differentiation.  What 
we see in Eqs. (\ref{b1}) and (\ref{b2}) is one statement of that 
relationship. But there is also the alternate form, which is stated 
for real variables as
\begin{equation}
\frac{d}{dt}\; \int^{t}\; f(t')\;dt' = f(t).\label{b9}
\end{equation}
For our quaternionic variables we start by looking at
\begin{equation}
{\cal{D}}_{x}\; \int ^{x}\; {\cal{D}}_{x'}\;F(x') \label{b10}
\end{equation}
and then apply the first differential operator to the coordinate $x$
in two parts: first the $\delta_{\parallel}$ part and then the 
$\delta_{\perp}$ part. The result is just the integrand evaluated at 
the point $x$:
\begin{equation}
 = {\cal{D}}_{x}\;F(x); \label{b11}
\end{equation}
and this is just what we should expect from the right hand side of 
Eq. (\ref{b2}), with $x_{b}$ replaced by $x$.

\section{Discussion}%3

Following what was stated in the earlier work, \cite{1}, we do require the functions 
$F(x)$ to be real analytic functions along the path of integration.

Our main result Eq. (\ref{b2}) implies that the result of the 
integration depends only on the end points and is independent of the 
path. This is true if we also require that the function $F(x)$ be 
single valued. Then, we have the result that the integral over any 
closed path, ending up at the same point where it started, is zero. 
For more discussion of this, see Appendix B.

\vskip 0.5cm
\setcounter{equation}{0}
\def\theequation{A.\arabic{equation}}
\boldmath
\noindent{\bf Appendix A}%A
\unboldmath
\vskip 0.5cm

If we look at the real integral and try to guess how 
to generalize it to the non-commutative quaternions, we might start with, 
\begin{equation}
\int f(t)\;dt  \; \stackrel{?}{\longrightarrow} \;  \frac{1}{2}\;\int\; (dx\; F(x) + 
F(x)\;dx) \;;\label{A1}
\end{equation}
but why should $dx$ only appear on the outside; why not also in the 
middle of the function $F(x)$?

Let's try a most symmetrical arrangement with the function $F(x) = 
x^{n}$:
\begin{equation}
\int t^{n}\; dt \;\stackrel{?}{\longrightarrow} \; 
\frac{1}{n+1}\;\int\;(dx\; x^{n} + 
x\;dx\;x^{n-1} + x^{2}\;dx\;x^{n-2} + \ldots + x^{n}\;dx).\label{A2}
\end{equation}
But we can recognize that the long expression in parentheses on the 
right hand side of this is nothing other than ${\cal{D}} x^{n+1}$:
\begin{equation}
{\cal{D}} F(x) \equiv  F(x+dx) -F(x), \;\; to\; first\; order\; in \; 
dx .\label{A2.1}
\end{equation}
So we would then write,
\begin{equation}
\int t^{n}\; dt\; \longrightarrow\; \frac{1}{n+1}\int {\cal{D}} 
x^{n+1} = \frac{x^{n+1}}{n+1},\label{A3}
\end{equation}
using our defining Eq. (\ref{b2}). Now, this looks quite familiar.

We can extend this to any power series and thus offer the following 
rule.  For any analytic function of a real variable $f(t)$, for which 
we know the integral,
\begin{equation}
\int\;f(t)\;dt = h(t),\label{A4}
\end{equation}
we can make the correspondence to quaternionic integration as follows:
\begin{equation}
\int\;f(t)\;dt \;\longrightarrow \;\int\;{\cal{D}}h(x) = 
h(x).\label{A5}
\end{equation}
While this may look trivial for real and complex variables, it is 
something new for non-commuting variables. This is because we have 
carefully defined and   
studied  the operator ${\cal{D}}$.

\vskip 0.5cm
\setcounter{equation}{0}
\def\theequation{B.\arabic{equation}}
\boldmath
\noindent{\bf Appendix B}%B
\unboldmath
\vskip 0.5cm

In the familiar study of functions of a complex variable, we have the 
 rule about integrals around a pole,
\begin{equation}
\oint \frac{g(z)}{z}\; dz = 2 \pi i\; g(0) \label{B1}
\end{equation}
for well behaved functions $g$. How does this square with our general 
statement above that any integral over a closed path would be zero?

Our formulation above, as applied to a complex variable, reads
\begin{equation}
\int_{a}^{b}\; \frac{df(z)}{dz}\; dz = f(b) - f(a)\label{B2}
\end{equation}
and this should be zero when $f(z)$ is analytic and single valued over 
the contour of integration and that contour is closed: $a=b$. The way 
to get the integral (\ref{B2}) to look like the integral (\ref{B1}) 
is to choose the function $f(z) = ln\; z$. But this function is not 
single valued. In going once around the origin, $ln\; z$ changes by 
exactly the amount $2 \pi i$.  So we have agreement between the 
results of integrals (\ref{B1}) and (\ref{B2}).

\end{document}